\documentclass[a4paper,draft]{amsart}
\usepackage{amsmath, amssymb,braket}
\usepackage{amsthm}
\usepackage[vcentermath,enableskew]{youngtab}
\usepackage{bm}
\newcommand{\partitionof}{\vdash}
\newcommand{\tabloid}[1]{\boldsymbol\{ #1 \boldsymbol\}}
\newcommand{\Tabloidson}[1]{\mathbb{T}_{#1}}

\newcommand{\TT}{\boldsymbol{T}}
\newcommand{\Ss}{\boldsymbol{S}}
\newcommand{\Tt}{\boldsymbol{t}}
\newcommand{\Ll}{\boldsymbol{l}}
\newcommand{\Mm}{\boldsymbol{m}}

\newcommand{\marked}{\ast}
\newcommand{\boldI}{{1^{\marked}}}
\newcommand{\boldII}{{2^{\marked}}}
\newcommand{\boldIII}{{3^{\marked}}}
\newcommand{\boldIV}{{4^{\marked}}}
\newcommand{\boldV}{{5^{\marked}}}

\newcommand{\ggamma}{\boldsymbol{\gamma}}
\newcommand{\mmu}{\boldsymbol{\mu}}

\newcommand{\ZZ}{\mathbb{Z}}
\newcommand{\NN}{\mathbb{N}}
\newcommand{\CC}{\mathbb{C}}
\newcommand{\numof}[1]{\left| #1 \right|}

\newcommand{\induce}[3]{\operatorname{Ind}_{#1}^{#2}{#3}}
\newcommand{\charof}{\operatorname{Char}}

\newcommand{\defit}[1]{{\em #1}}

\theoremstyle{plain}
\newtheorem{thm}{Theorem}[section]
\newtheorem{lemma}[thm]{Lemma}
\newtheorem{cor}[thm]{Corollary}
\newtheorem{theorem}[thm]{Theorem}

\theoremstyle{remark}
\newtheorem{remark}[thm]{Remark}

\newtheorem{example}[thm]{Example}

\theoremstyle{definition}
\newtheorem{definition}[thm]{Definition}

\allowdisplaybreaks[3]

\begin{document}
\title[Tabloids and Weighted Sums of Characters]{
Tabloids and Weighted Sums of \\
 Characters of Certain Modules of \\
the Symmetric Groups}
\author[Numata, y.]{Numata, yasuhide}
\thanks{Department of Mathematics, Hokkaido University. e-mail:{\tt nu@math.sci.hokudai.ac.jp}}
\email{nu@math.sci.hokudai.ac.jp}
\address{
  Department of Mathematics\\
  Hokkaido University\\
  Sapporo 060-0810\\
  Japan}

\begin{abstract}
We consider certain modules of the symmetric groups whose basis 
elements are called tabloids. Some of these modules are isomorphic 
to subspaces of the cohomology rings of subvarieties of 
flag varieties as modules of the symmetric groups. We give a
 combinatorial description for some weighted sums of their characters,
 i.e., we introduce 
combinatorial objects called $(\rho,\Ll)$-tabloids and rewrite 
 weighted sums of 
characters as the numbers of these combinatorial objects. 
We also 
consider the meaning of these combinatorial objects, i.e., we construct 
a correspondence between  $(\rho,\Ll)$-tabloids and tabloids whose images are eigenvectors 
of the action of an element of cycle type $\rho$ 
in quotient modules. 
\end{abstract}

\maketitle

\section{Introduction}

Let $W$ be a finite reflection group.
In some $\ZZ$-graded $W$-modules $R=\bigoplus_{d} R^{d}$, we have 
a phenomenon called ``coincidence of dimensions'' 
(\cite{m}, \cite{MN} and so on), i.e.,
some integers $l$  satisfy
the equations
\begin{gather*}
\dim \bigoplus_{i\in\ZZ} R^{i l+k} = \dim \bigoplus_{i\in\ZZ} R^{i l+k'}
\end{gather*}
for all $k$ and $k'$.
Induced modules often give a proof of the phenomenon.
More specifically,
let a subgroup $H(l)$ of $W$ and $H(l)$-modules $z(k;l)$ satisfy
\begin{align*}
\bigoplus_{i\in\ZZ} R^{i l+k} \simeq \induce{H(l)}{W}{z(k;l)},  && \dim z(k;l)=\dim z(k';l)
\end{align*}
for all $k$ and $k'$, where $ \induce{H(l)}{W}{z(k;l)}$ denotes the
induced module.
Since 
\begin{align*}
\dim \induce{H(l)}{W}{z(k;l)} = \numof{W/H(l)} \cdot \dim z(k;l),
\end{align*}
we can prove the phenomenon
by the datum $(H(l),\Set{z(k;l)})$.

We consider the case where  $W$ is the $m$-th symmetric group $S_m$
and $R$ are the $S_m$-modules $R_\mu$ called Springer modules. 
The Springer modules $R_\mu$ are graded algebras parametrized by
 partitions $\mu \partitionof m$.
As $S_m$-modules, $R_\mu$ are  
isomorphic to cohomology rings
of the variety of the flags
 fixed by a unipotent matrix the sizes of which Jordan blocks are $\mu$.
  (See \cite{HotSpring, Sp76, Sp78}. See also \cite{dp, tani} for
algebraic construction.)
Let $\mu$ be an $l$-partition,
where an $l$-partition means a partition whose multiplicities are divisible
by $l$.
To prove 
coincidence of dimensions of the Springer module $R_\mu$ 
by induced modules,
Morita-Nakajima \cite{mn} explicitly calculated 
the Green polynomials corresponding to $\mu$
at $l$-th roots of unity.
These values are nonnegative integers.

Our first motivation for this paper is to describe these nonnegative values of 
the Green polynomials as  numbers of some combinatorial objects.
Our second motivation is to give a meaning of the
combinatorial objects
in terms of modules 
 $\induce{H_\mu(l)}{S_m}{Z_\mu(k;l)}$ in Morita-Nakajima \cite{mn}.
For these purposes, 
 we introduce some $S_m$-modules,
which are  realizations of  $\induce{H_\mu(l)}{S_m}{Z_\mu(k;l)}$ for
special parameters, and give a combinatorial description for  weighted
sums of  their characters.

In Section $\ref{defsec}$, we introduce $S_m$-modules
$M^{\mmu}$ and their quotient modules $M^{\mmu}(k;\Ll)$
for some $n$-tuple $\mmu$ of Young diagrams.
When $n=1$,
this module $M^{(\mu)}(k;(l))$ is a realization of $\induce{H_\mu(l)}{S_m}{Z_\mu(k;l)}$
 in \cite{mn}.
We also introduce combinatorial objects called marked $(\rho,\Ll)$-tabloids
to describe weighted sums of characters of  $M^{\mmu}(k;\Ll)$.
When $n=1$, the number of marked $(\rho,(l))$-tabloids 
coincides with the right hand side of 
the explicit formula $(3.1)$ of Green polynomials  in \cite{mn}.
Our main result is the description of a weighted sum 
\begin{gather*}
\sum_{k \in \ZZ/ l \ZZ} \zeta_l^{j k} \charof \left(M^{\mmu}(k;\Ll) \right)(\sigma)
\end{gather*}
of characters of
$M^{\mmu}(k;\Ll)$ as the number of marked $(\rho,\ggamma)$-tabloids on
$\mmu$ 
for the primitive  $l$-th root $\zeta_l$ of unity and $\sigma\in S_m$
of cycle type $\rho$
in Section $\ref{mainsec}$.
We prove the main result 
in Section $\ref{proofsec}$
by constructing bijections.

\section{Notation and Definition}\label{defsec}

We identify a partition $\mu=(\mu_1\geq\mu_2\geq\cdots)$ of $m$ with its Young diagram 
$\Set{(i,j) \in \NN^2 |1 \leq j \leq \mu_i}$ with $m$ boxes.
If $\mu$ is a Young diagram with
$m$ boxes, 
we write $\mu \partitionof m$ and identify a Young diagram $\mu$ with 
the array of $m$ boxes having left-justified rows with the $i$-th row
containing  $\mu_i$ boxes; for example,
\begin{gather*}
(2,2,1) = \yng(2,2,1)  \partitionof 5.
\end{gather*}
For an integer $l$, a Young diagram $\mu$ is called an \defit{$l$-partition}
if multiplicities $m_i=\numof{\Set{k|\mu_k=i}}$ of $i$
are divisible by $l$ for all $i$.

Let $\mu$ be a Young diagram with $m$ boxes. 
We call a map  $T$ a \defit{numbering on $\mu$ with $\{1,\ldots,n\}$} 
if $T$ is an injection $\mu \ni (i,j) \mapsto T_{i,j} \in \Set{1,\ldots , n}$.
We identify a map $T:\mu\to\NN$
with a diagram putting $T_{i,j}$  in
each box in the $(i,j)$ position; for example, 
\begin{gather*}
\young(23,61,5)
\end{gather*}
is identified with a  numbering $T$ on $(2,2,1)$ which maps $T_{1,1}=2$,
$T_{1,2}=3$, $T_{2,1}=6$, $T_{2,2}=1$, $T_{3,1}=5$.
For $\mu \partitionof m$,
 $t_\mu$ denotes the numbering which maps
$(t_\mu)_{i,j}=j+\sum_{k=0}^{i-1}\mu_k$; i.e.,
the numbering obtained by putting numbers from $1$ to $m$ on boxes of $\mu$
from left to right in each row, starting in the top row and moving
to the bottom row.
For example, 
\begin{gather*}
t_{\ \tiny\yng(2,2,1)}=\young(12,34,5).
\end{gather*}

Two numberings $T$ and $T'$  on $\mu \partitionof m$ are said to be
\defit{row-equivalent} if their corresponding rows
consist of the
same elements.
We call a row-equivalence class $\tabloid{ T }$ a \defit{tabloid}.

Let $T$ be a numbering  $T$  on a Young diagram 
$\mu \partitionof m$ with $\{1,\ldots,n\}$.
Then $\sigma \in S_n$ acts on  $T$  from left
as $(\sigma T)_{i,j}=\sigma (T_{i,j})$.
For example,
\begin{gather*}
(1,2,3,4)\ \young(12,3)=\young(23,4).
\end{gather*}
This left action induces a left action on tabloids by 
$\sigma \tabloid{ T } = \tabloid{ \sigma T }$.

For a numbering $T$ on  $\mu \partitionof m$ with $\{1,\ldots,n\}$,
we define $S_T$ to be the subgroup 
\begin{gather*}
S_{\Set{T_{1,1},T_{1,2},\ldots,T_{1,\mu_1}}} \times
 S_{\Set{T_{2,1},T_{2,2},\ldots,T_{2,\mu_2}}} \times \cdots 
\end{gather*}
of the $n$-th symmetric group $S_n$,
where $S_{\Set{i_1,\ldots,i_k}}$ denotes the symmetric group of the
letters $\Set{i_1,\ldots,i_k}$.
It is obvious that
$S_T$ and the Young subgroup $S_\mu$ are isomorphic as
groups
for a numbering $T$ on $\mu \partitionof m$.
It is also clear that
$\sigma \tabloid{ T } = \tabloid{  T }$
for $\sigma\in S_T$.

For a numbering $T$ on an $l$-partition  $\mu \partitionof m$,
we define $a_{T,l}$ to be the product 
\begin{gather*}
\prod_{(l i+1,j)\in\mu}(T_{l i+1,j}, T_{l i+2,j},\ldots, T_{l i+l,j})
\end{gather*}
of $m/l$ cyclic permutations of
length $l$.
For example, $a_{t_{(2,2,1,1)},2}=(13)(24)(56)$.
We write $a_{\mu,l}$ for $a_{t_{\mu},l}$.

Let $\mu$ be an $l$-partition of $m$ and 
$\Braket{a_{\mu,l}}$ the cyclic group of order $l$ generated by
$a_{\mu,l}$.
For each  numbering  $T$  on $\mu$ with $\{1,\ldots,n\}$,
there exists $\tau_T \in S_{n}$ such that
$T=\tau_T t_\mu$.
Since the map $\tau_{T}|_{\Set{1,\ldots,m}}$ restricting $\tau_{T}$
to $\Set{1,\ldots,m}$ is unique,
 $\sigma \in \Braket{a_{\mu,l}}$ acts on   
$T$   from right
as $T \sigma = \tau_T \sigma t_\mu$.
For each numbering $T$ on an $l$-partition $\mu$,
the $(\overline{r}+l q)$-th row of $T a_{\mu,l}$
is  the $(\overline{r+1}+l q)$-th row of  $T$,
where $\overline{r}$ and $\overline{r+1}$ are in $\ZZ/l\ZZ=\Set{1,\ldots,l}$.
This right action also induces a right action on tabloids by 
$ \tabloid{ T } \sigma  = \tabloid{ T  \sigma }$.

In this paper, we consider $n$-tuples of Young diagrams.
Throughout this paper,
let  $\Mm=(m_{1},m_{2},\ldots,m_{n})$ and  $\Ll=(l_{1},l_{2},\ldots,l_{n})$
be $n$-tuples of positive integers,
$m$ the sum $\sum_h m_h$, 
 $l$  the least common multiple of $\Set{l_{i}}$
and 
$\zeta_{k}$ the primitive $k$-th root of unity.
We call an $n$-tuple $\mmu=(\mu^{(1)},\mu^{(2)},\ldots,\mu^{(n)})$ of Young
diagrams an \defit{$\Ll$-partition of $\Mm$}
if $\mu^{(h)}$ is an $l_{h}$-partition of $m_{h}$ for each $h$.
We identify an $\Ll$-partition $\mmu$ with
the disjoint union 
$\coprod_{h} \mu^{(h)} =\Set{(i,j;h)|(i,j) \in \mu^{(h)}}$
of Young diagrams $\mu^{(h)}$.
We call an $n$-tuple $\TT=(T^{(1)},\ldots,T^{(n)})$ of numberings $T^{(h)}$ on $\mu^{(h)}$
\defit{a numbering on an $\Ll$-partition $\mmu$} 
if the map  $\TT:\mmu \ni (i,j;h) \mapsto T_{i,j}^{(h)} \in\Set{1,\ldots, m}$ 
is bijective.
For an $\Ll$-partition $\mmu$ of $\Mm$,
$\Tt_{\mmu}$ denotes the $n$-tuple of the numberings $t_{\mmu}^{(h)}$ 
which maps 
$(i,j;h)$ to 
$(t^{(h)}_{\mmu})_{i,j}=(t_{\mu^{(h)}})_{i,j}+\sum_{k=1}^{h-1} m_k$; 
i.e.,
$\Tt_{\mmu}$ is
the numbering on an $\Ll$-partition $\mmu$
obtained by putting numbers from $1$ to $m$ on boxes of $\mmu$
from left to right in each row, starting in the top row and moving
to the bottom in each Young diagram, starting from $\mu^{(1)}$ to $\mu^{(n)}$.
For example, 
\begin{gather*}
\Tt_{\tiny \left(\ \yng(2,2,1,1)\ ,\ \yng(2,1)\ \right)}=\left(\
 \young(12,34,5,6)\ ,\ \young(78,9)\ \right).
\end{gather*}

Let $\TT$  be a numbering on an $\Ll$-partition $\mmu$ of $\Mm$.
We define $S_{\TT}$ to be the subgroup 
$S_{\TT} = S_{T^{(1)}} \times S_{T^{(2)}} \times\cdots\times S_{T^{(n)}}$
of $S_{m}$.
The subgroup  $S_{\TT}$ and
the Young subgroup
$S_{\overline{\mmu}}$  are isomorphic as groups,
where $\overline{\mmu}$ is the partition  obtained from
$(\mu^{(1)}_{1},\mu^{(2)}_{1},\ldots,\mu^{(n)}_{1},\mu^{(1)}_{2},\mu^{(2)}_{2},\ldots,\mu^{(n)}_{2},\ldots)$
by sorting  in descending order.
We write $S_{\mmu}$ for $S_{\Tt_{\mmu}}$.
We define $a_{\TT,\Ll}$ to be $a_{T^{(1)},l_{1}}\cdot a_{T^{(2)},l_{2}}\cdots a_{T^{(n)},l_{n}}$.
We write $a_{\mmu,\Ll}$ for $a_{\Tt_{\mmu},\Ll}$

Two numberings $\TT$ and $\Ss$ 
on an $\Ll$-partition of $\Mm$
are said to be \defit{row-equivalent} 
if $T^{(h)}$ and $S^{(h)}$ are row-equivalent for each $h$.
The set of numberings 
whose  components are arranged in ascending order 
in each row is a complete set of representatives for row-equivalence classes.
A row-equivalence class  of a numbering
$\TT$ on an $\Ll$-partition $\mmu$
is an $n$-tuple 
$(\tabloid{T^{(1)}},\tabloid{T^{(2)}},\ldots,\tabloid{T^{(n)}})$ of
tabloids $\tabloid{T^{(h)}}$ on $\mu^{(h)}$.
We also call a row-equivalence class 
$(\tabloid{T^{(h)}})$ of numbering
$\TT$ on an $\Ll$-partition
a \defit{tabloid on an $\Ll$-partition}.
We also write $\tabloid{\TT}$ for  $(\tabloid{T^{(h)}})$.
The set of all tabloids on an $\Ll$-partition $\mmu$
of $\Mm$ is
denoted by $\Tabloidson{\mmu}$.
We define $M^{\mmu}$ to be the $\CC$-vector space 
$\CC \Tabloidson{ \mmu }$ whose basis is the set $\Tabloidson{ \mmu }$ of tabloids on 
$\mmu$.

Let $\TT$ be a numbering 
on an $\Ll$-partition of $\Mm$.
Then $\sigma \in S_{m}$ acts on  $\TT$
from left as $\sigma (T^{(h)})=(\sigma T^{(h)})$.
This left action induces a left action on tabloids by 
$\sigma \tabloid{\TT} = \tabloid{\sigma \TT}$.
For the partition $\overline{\mmu}\partitionof m$, 
$M^{\mmu}$ and $\induce{S_{\overline{\mmu}}}{S_m}{1}$ are isomorphic as
left $S_m$-modules,
where  $1$ denotes the 
 the trivial module of the Young subgroup $S_{\overline{\mmu}}$. 

Let $\TT$ be a numbering on  an $\Ll$-partition $\mmu$ of $\Mm$.
Since there uniquely exists $\tau_{\TT} \in S_{m}$ such that
$\TT=\tau_{\TT} \Tt_{\mmu}$,
$\sigma \in \Braket{a_{\mmu,\Ll}}$ acts on $\TT$  
from right as $\TT \sigma = \tau_{\TT} \sigma \Tt_{\mmu}$.
This right action also induces  a right action on tabloids by 
$ \tabloid{ \TT } \sigma  = \tabloid{ \TT  \sigma }$.

Next we introduce $S_m$-modules $M^{\mmu}(k;\Ll)$,
one of main objects in this paper.
We need some definitions to introduce $M^{\mmu}(k;\Ll)$.
\begin{definition}

Let 
$\Tabloidson{\mmu}^{\Ll}$ be  the  subset
$\Set{ a_{\mmu,\Ll}^i \tabloid{\Tt_{\mmu}} | i \in \ZZ / l \ZZ}$
of tabloids
for an $\Ll$-partition $\mmu$ of $\Mm$.
We define $Z_{\mmu}(\Ll)$ to be the $\CC$-vector space 
$\CC\Tabloidson{\mmu}^{\Ll}$ whose basis is 
$\Tabloidson{\mmu}^{\Ll}$.
This $l$-dimensional vector space 
 is 
a left module of the semi-direct product 
$S_{\mmu} \rtimes \Braket{ a_{\mmu,\Ll}}$ and 
a right module of the cyclic group 
$\Braket{a_{\mmu,\Ll}}$ of order $l$.

For $k\in \ZZ/l\ZZ$, let $I_{\mmu}(k;\Ll)$ denote
  the submodule of  $Z_{\mmu}(\Ll)$  generated by 
\begin{gather*}
\Set{ a_{\mmu,\Ll}^i \tabloid{\Tt_{\mmu}} -
 \zeta_l^{k i} \tabloid{\Tt_{\mmu}} | i \in \ZZ/l\ZZ}.
\end{gather*}
We define  $Z_{\mmu}(k;\Ll)$ to be the quotient module
\begin{gather*}
Z_{\mmu}(\Ll)/I_{\mmu}(k;\Ll).
\end{gather*}
For each $k$, $Z_{\mmu}(k;\Ll)$ is a  one-dimensional left module of
the semi-direct product 
$S_{\mmu} \rtimes \Braket{ a_{\mmu,\Ll}}$.
This 
left $S_{\mmu} \rtimes \Braket{ a_{\mmu,\Ll}}$-module 
$Z_{\mmu}(k;\Ll)$ is generated by $\tabloid{\Tt_{\mmu}}$,
and $a_{\mmu,\Ll}$ acts
on $\tabloid{\Tt_{\mmu}}$
by 
\begin{gather*}
 a_{\mmu,\Ll} \tabloid{\Tt_{\mmu}} 
=\zeta_{l}^k \tabloid{\Tt_{\mmu}}
\end{gather*}
in $Z_{\mmu}(\Ll)/I_{\mmu}(k;\Ll)$.

Let $\widetilde I_{\mmu}(k;\Ll)$ be $\CC[S_{m}] I_{\mmu}(k;\Ll)$. 
Finally, we define an $S_{m}$-module  $M^{\mmu}(k;\Ll)$ 
to be 
\begin{gather*}
M^{\mmu}/\widetilde I_{\mmu}(k;\Ll).
\end{gather*}
By definition, 
the $S_n$-module $M^{\mmu}(k;\Ll)$ is a realization of the induced module 
$\induce{S_{\mmu} \rtimes \Braket{ a_{\mmu,\Ll}}}{S_m}  Z_{\mmu}(k;\Ll)$.
\end{definition}

\begin{remark}
For an $l$-partition $\mu$ of $m$, 
our module $M^{(\mu)}(k;(l))$ gives a realization of
the $S_m$-module  $\induce{H_\mu(l)}{S_m}{Z_\mu(k;l)}$ in Morita-Nakajima \cite{mn}.
For $n$-tuple $\Set{l_h}$ of integers, $M^{\mmu}(k;\Ll)$ is a realization of
the induced module 
\begin{gather*}
\induce{S_{m_1}\times\cdots\times S_{m_n}}{S_{m_1+\cdots +m_n}}{M^{\mu^{(1)}}(k;l_1) \otimes \cdots \otimes M^{\mu^{(n)}}(k;l_n)},
\end{gather*}
where $M^{\mu}(k;l)$ denotes $M^{(\mu)}(k;(l))$.
\end{remark}

\begin{remark}
Since 
$\widetilde I_{\mmu}(k;\Ll) = \CC[S_{m}] I_{\mmu}(k;\Ll)$
is generated by
\begin{gather*}
\Set{ \tau a_{\mmu,\Ll}^i \tabloid{\Tt_{\mmu}} -\zeta_{l}^{i k} \tau \tabloid{\Tt_{\mmu}}| i \in \ZZ / l\ZZ , \tau \in S_{m}},
\end{gather*}
$\widetilde I_{\mmu}(k;\Ll)$ is also generated by 
\begin{gather*}
\Set{ \tabloid{\TT} a_{\mmu,\Ll}^i  -\zeta_{l}^{i k}\tabloid{\TT} | \tabloid{\TT} \in \Tabloidson{\mmu} , i\in \ZZ/l\ZZ}.
\end{gather*}
Hence $a_{\mmu,\Ll}$ acts on tabloids $\tabloid{\TT}$ by 
\begin{gather*}
\tabloid{\TT}a_{\mmu,\Ll}=\zeta_{l}^{k}\tabloid{\TT}
\end{gather*}
in $M^{\mmu}(k;\Ll)$.
\end{remark}

 We introduce the following combinatorial objects to describe the
characters of $M^{\mmu}(k;\Ll)$. 
\begin{definition}
For a Young diagram $\rho \partitionof m$,
we call a map $Y:\mmu \to \NN$
a \defit{$(\rho,\Ll)$-tabloid} on an $\Ll$-partition $\mmu$ of $\Mm$
if the following are satisfied:
\begin{itemize}
\item $\numof{Y^{-1}(\Set{k})} = \rho_k$ for all $k$,
\item for each $k$, there exist $h\in\NN$ and  $(i',j') \in \NN^2$
such that $\rho_k$ is divisible by $l_{h}$ and
\begin{gather*}
Y^{-1}(\Set{k})=\Set{ (i+i', j+j';h)| (i,j) \in \left(\left(\frac{\rho_k}{l_{h}}\right)^{l_{h}}\right)
      \partitionof \rho_k},
\end{gather*}
\item for each $(i,j;h)$, $(i,k;h) \in \mmu$, $Y(i,j;h) \leq Y(i,k;h)$ if $j \leq k$.
\end{itemize}
\end{definition}

\begin{example}
For example, 
\begin{gather*}
\left(\ \young(34,34,11,11)\ ,\ \young(225)\ \right)
\end{gather*}
is a $((4,2,2,2,1),(2,1))$-tabloid on
$((2,2,2,2),(3))$.
\end{example}

\begin{definition}
We call a pair $(Y,c)$ 
a \defit{marked $(\rho,\Ll)$-tabloid}  on an $\Ll$-partition $\mmu$ of $\Mm$
if the following are satisfied:
\begin{itemize}
\item $Y$ is a $(\rho,\Ll)$-tabloid  on an $\Ll$-partition $\mmu$,
\item $c$ is a map from $\Set{i|\rho_i\neq 0}$ to $\coprod_{h} \ZZ/l_{h}\ZZ$,
\item $c(i)$ is in $\ZZ/l_{h}\ZZ$ if  $Y^{-1}(\Set{i}) \subset \mu^{(h)}$.  
\end{itemize}
\end{definition}

For a marked $(\rho,\Ll)$-tabloid $(Y,c)$, the inverse image $Y^{-1}(\Set{i})$
has $l_h$ rows 
and $c({i})$ is in $\ZZ/l_h\ZZ=\Set{1,\ldots,l_h}$
if $Y^{-1}(\Set{i})$ is in $\mu^{(h)}$.
We identify $(Y,c)$ with the diagram obtained from the diagram of $Y$
by putting $\marked$ in the left-most box of 
the $c(i)$-th row of the inverse image $Y^{-1}(\Set{i})$,
where we identify $\ZZ/l_h\ZZ$ with the set $\Set{1,\ldots,l_h}$ of complete representatives.

\begin{example}
Let 
\begin{gather*}
Y=\left(\ \young(34,34,11,11) \ ,\ \young(225) \ \right)
\end{gather*}
and let $c$ be the map such that $c(1)=2$,
$c(3)=1$, $c(4)=2 \in \ZZ/2\ZZ$ and
$c(2)=c(5)=1 \in \ZZ/ 1 \ZZ$,
then $\left(Y,c\right)$
is a marked $(2,1)$-tabloid.
We write  
\begin{gather*}
\left(\ \young(\boldIII 4,3\boldIV,11,\boldI 1)\ ,\ \young(\boldII
 2\boldV)\ \right)
\end{gather*}
for $(Y,c)$.
\end{example}

\begin{remark}
It follows from a direct calculation that
the number of marked $(\rho,(l))$-tabloids on an $(l)$-partition $(\mu)$
equals the right hand side of the equation $(3.1)$ in Morita-Nakajima
\cite{mn}. 
\end{remark}

\begin{definition}
Let 
 $\mmu$ be an $\Ll$-partition of $\Mm$
and   $\ggamma=(\gamma_h)$  an $n$-tuple of 
 integers such that $l_{h}$ is divisible by $\gamma_{h}$.
For a Young diagram $\rho \partitionof m$,
we call a map $Y:\mmu \to \NN$
a \defit{$(\rho,\ggamma,\Ll)$-tabloid} on  $\mmu$
if the following are satisfied:
\begin{itemize}
\item $\numof{Y^{-1}(\Set{k})} = \rho_k$ for all $k$,
\item for each $k$, there exist  $h$ and $(i',j') \in \NN^2$
such that $\rho_k$ is divisible by $\gamma_{h}$ and
\begin{gather*}
Y^{-1}(\Set{k})=\Set{ \left(\frac{i l_{h}}{\gamma_{h}}+i', j+j';h\right) | (i,j) \in \left(\left(\frac{\rho_k}{\gamma_{h}}\right)^{\gamma_{h}}\right)
      \partitionof \rho_k},
\end{gather*}
\item for each $(i,j;h)$, $(i,k;h) \in \mmu$, $Y(i,j;h) \leq Y(i,k;h)$ if $j \leq k$.
\end{itemize}
\end{definition}

\begin{example}
For example, 
\begin{gather*}
\left(\ \young(2224,1111,2224,1111)\ ,\ \young(33,55,6)\ \right)
\end{gather*}
is an $((8,6,2,2,2,1),(2,1),(4,1))$-tabloid on
$((4,4,4,4),(5))$.
\end{example}

A $(\rho,\Ll,\Ll)$-tabloid on an $\Ll$-partition $\mmu$
is a $(\rho,\Ll)$-tabloid on $\mmu$.

For an $\Ll$-partition $\mmu$ and an $n$-tuple $\ggamma=(\gamma_h)$
such that $l_h$ is divisible by $\gamma_h$,
it follows that
\begin{gather*}
\numof{\Set{Y|\text{a $(\rho,\ggamma,\Ll)$-tabloid on $\mmu$}}}
=\numof{\Set{Y|\text{a $(\rho,\ggamma)$-tabloid on $\mmu$}}}. 
\end{gather*}

\begin{definition}
We call a pair $(Y,c)$ 
a \defit{marked $(\rho,\ggamma,\Ll)$-tabloid} on an $\Ll$-partition
 $\mmu$ of $\Mm$ 
 the following are satisfied:
\begin{itemize}
\item $Y$ is a $(\rho,\ggamma,\Ll)$-tabloid  on an $\Ll$-partition $\mmu$,
\item $c$ is a map from $\Set{i|\rho_i\neq 0}$ to $\coprod_{h}\ZZ/\gamma_{h}\ZZ$,
\item $c(i)$ is in $\ZZ/\gamma_{h}\ZZ$ if  $Y^{-1}(\Set{i}) \subset \mu^{(h)}$.  
\end{itemize}
\end{definition}

Similar to the case of  marked $(\rho,\Ll)$-tabloids, we identify  a marked $(\rho,\ggamma,\Ll)$-tabloid $(Y,c)$ with the diagram obtained from the diagram of $Y$
by putting $\marked$ in the left-most box of 
the $c(i)$-th row of the inverse image $Y^{-1}(\Set{i})$.

\begin{example}
For example, 
\begin{gather*}
\left(\ \young(\boldIII 4,11,3\boldIV,\boldI 1)\ ,\ \young(\boldII 2\boldV)\ \right)
\end{gather*}
is a marked $((4,2,2,2,1),(2,1),(4,1))$-tabloid.
\end{example}

For an $\Ll$-partition $\mmu$ and an $n$-tuple $\ggamma=(\gamma_h)$
such that $l_h$ is divisible by $\gamma_h$,
it follows that
\begin{gather}
\begin{split}
&\numof{\Set{(Y,c)|\text{a marked $(\rho,\ggamma,\Ll)$-tabloid on $\mmu$}}}\\
=&\numof{\Set{(Y,c)|\text{a marked $(\rho,\ggamma)$-tabloid on $\mmu$}}}.
\end{split}
\label{eqofnum}\end{gather}

\section{Main Results}\label{mainsec}
The following are the main results of this paper, 
proved in Section $\ref{proofsec}$.

\begin{theorem}\label{gencharathm}
For an integer $j$, let $\mmu$ be an $\Ll$-partition 
and $\ggamma$ an $n$-tuple of integers such that 
$\gamma_{h}$ is the order of $ \zeta_{l_h}^j $.
For $\sigma \in S_{m}$ of cycle type $\rho$,
\begin{gather*}
\sum_{k \in \ZZ/ l \ZZ} \zeta_l^{j k} \charof \left(M^{\mmu}(k;\Ll) \right)(\sigma)
=\numof{\Set{(Y,c)|\text{a marked $(\rho,\ggamma)$-tabloid on $\mmu$}}}.
\end{gather*}
\end{theorem}

\begin{theorem}\label{genidenthm}
For an integer $j$, let $\mmu$ be an $\Ll$-partition 
and $\ggamma$ an $n$-tuple of integers such that 
$\gamma_{h}$ is the order of $ a_{\mu^{(h)},l_{h}}^{j} $.
Tabloids $\TT$ on $\mmu$ satisfying 
$\sigma \tabloid{\TT} = \tabloid{\TT} a_{\mmu,\Ll}^{-j}$
are parameterized by 
marked $(\rho_\sigma,\ggamma,\Ll)$-tabloids on $\mmu$,
where $\rho_\sigma$ is the  cycle type of $\sigma$.
\end{theorem}

 Corollary $\ref{charathm}$ and Corollary $\ref{identhm}$ below
are obtained from Theorem  $\ref{gencharathm}$ and Theorem  $\ref{genidenthm}$
by applying $j=1$.

\begin{cor}\label{charathm}
For $\sigma \in S_{m}$ of cycle type $\rho$ and
an $\Ll$-partition $\mmu$, it follows that
\begin{gather*}
\sum_{k \in \ZZ/ l \ZZ} \zeta_l^k \charof \left(M^{\mmu}(k;\Ll) \right)(\sigma)
=\numof{\Set{(Y,c)|\text{a marked $(\rho,\Ll)$-tabloid on $\mmu$}}}.
\end{gather*}
\end{cor}

\begin{cor}\label{identhm}
Let $\mmu$ be an $\Ll$-partition.
Tabloids $\tabloid{\TT}$ on $\mmu$ satisfying 
$\sigma \tabloid{\TT} = \tabloid{\TT} a_{\mmu,\Ll}^{-1}$
are parameterized by 
$\Ll$-fillings on $(\rho_\sigma,\Ll)$-tabloids on $\mmu$,
where $\rho_\sigma$ is the  cycle type of $\sigma$.
\end{cor}

The following corollary 
directly follows from Theorem $\ref{gencharathm}$.
\begin{cor}\label{chareqcor}
For an integer $j$, let $\mmu$ be an $\Ll$-partition, 
 $\ggamma$ an $n$-tuple of integers such that 
$\gamma_{h}$ is the order of $ \zeta_{l_h}^j $. 
For $\sigma \in S_{m}$ of cycle type $\rho$,
\begin{align*}
\sum_{k \in \ZZ/ l \ZZ} \zeta_l^{j k} \charof \left(M^{\mmu}(k;\Ll) \right)(\sigma)
&=\sum_{k \in \ZZ/ \gamma \ZZ} \zeta_\gamma^{k} \charof \left(M^{\mmu}(k;\ggamma) \right)(\sigma)
\\
&=\numof{\Set{(Y,c)|\text{a marked $(\rho,\ggamma)$-tabloid on $\mmu$}}}.
\end{align*}
\end{cor}

\begin{example}\label{runningex}
Let $\mmu=((2,2),(4))$ and $\Ll=(2,1)$.
First we consider the case where $j=1$.
In this case,
all marked $((4,2,2),\Ll)$-tabloids on $\mmu$ are the following:
\begin{align*}
&\left(\ \young(\boldI 1,11)\ ,\ \young(\boldII 2\boldIII 3)\ \right),
 & & \left(\ \young(11,\boldI 1)\ ,\ \young(\boldII 2\boldIII 3)\ \right),\\
&\left(\ \young(\boldII\boldIII,23)\ ,\ \young(\boldI 111)\ \right),
 & & \left(\ \young(\boldII 3,2\boldIII)\ ,\ \young(\boldI 111)\ \right),\\
&\left(\ \young(2\boldIII,\boldII 3)\ ,\ \young(\boldI 111)\ \right),
 & & \left(\ \young(23,\boldII\boldIII)\ ,\ \young(\boldI 111)\ \right)
.\end{align*}
It follows from Corollary \ref{charathm} that
\begin{gather*}
\sum_{k \in \ZZ/ l \ZZ} \zeta_l^{k} \charof \left(M^{\mmu}(k;\Ll) \right)((1234)(56)(78))
=6.
\end{gather*}
Next consider the case where $j=2$.
Since $\zeta_{l_1}=\zeta_{2}=-1$ and $\zeta_{l_2}=\zeta_{1}=1$,
we have  $\gamma_1=\numof{\Braket{\zeta_{l_1}^2}} =1$ and
 $\gamma_2=\numof{\Braket{\zeta_{l_2}^2}} =1$.
All marked $((4,2,2),(1,1))$-tabloids on $\mmu$ are the following:
\begin{align*}
&\left(\ \young(\boldII 2,\boldIII 3)\ ,\ \young(\boldI 111)\ \right),
 & &\left(\ \young(\boldIII 3,\boldII 2)\ ,\ \young(\boldI 111)\ \right)
.\end{align*}
It follows from Theorem \ref{gencharathm} that
\begin{gather*}
\sum_{k \in \ZZ/ l \ZZ} \zeta_l^{2k} \charof \left(M^{\mmu}(k;\Ll) \right)((1234)(56)(78))
=2.
\end{gather*}
\end{example}

\section{Proof of Main Results}\label{proofsec}

In this section, we prove Theorem \ref{gencharathm} and Theorem \ref{genidenthm}.
First we prove that Theorem \ref{gencharathm} and Theorem
\ref{genidenthm} are equivalent; to prove the equivalence,
we prepare Lemma \ref{equivlemma}. Next we prove  Theorem \ref{genidenthm}
by giving an explicit parametrization; the  correspondence $\varphi$ defined in
Definition \ref{defofphi2} provides an explicit
parametrization.
We prove  Theorem \ref{genidenthm} first for the special element $\sigma_\rho$ of the cycle type
$\rho$, which is  Lemma \ref{correspondlemma2}.
We prepare Lemma \ref{eigenlemma} and Lemma \ref{shapelemma2}
to prove  Lemma \ref{correspondlemma2}.
Last, in Theorem \ref{lastlemma}, we generalize  Lemma
\ref{correspondlemma2}
for general elements of the cycle type $\rho$.
Theorem
\ref{lastlemma} is a realization of  Theorem \ref{genidenthm}.

First we prove the equivalence
of  Theorem \ref{gencharathm} and Theorem \ref{genidenthm}.
\begin{lemma}\label{equivlemma}
For an $\Ll$-partition $\mmu$ and $\sigma \in S_m$,
\begin{gather*}
\sum_{k \in \ZZ/ l \ZZ} \zeta_l^{k j} \charof \left(M^{\mmu}(k;\Ll) \right)(\sigma)
=\numof{\Set{ \tabloid{\TT} \in \Tabloidson{\mmu}|\sigma \tabloid{\TT} = \tabloid{\TT} a_{\mmu,\Ll}^{-j} }}.
\end{gather*}
\end{lemma}
\begin{proof}
Let $\overline{\Tabloidson{\mmu}}$ be a subset of tabloids
 $\Tabloidson{\mmu}$ whose image is a basis of
 $M^{\mmu}(0;\Ll)$.
The image of  $\overline{\Tabloidson{\mmu}}$ in $M^{\mmu}(k;\Ll)$
 is also a basis of   $M^{\mmu}(k;\Ll)$
for every $k$.
Let  $\Tabloidson{\mmu}^{(\sigma,j;\Ll)}$ be a subset 
\begin{gather*}
\Set{  \tabloid{\TT}\in \Tabloidson{\mmu} | \sigma \tabloid{\TT} = \tabloid{\TT} a^j_{\mmu,\Ll}}
\end{gather*} 
of $\Tabloidson{\mmu}$.
We write  $\overline{\Tabloidson{\mmu}^{(\sigma,j;\Ll)}}$ for
$\overline{\Tabloidson{\mmu}} \cap \Tabloidson{\mmu}^{(\sigma,j;\Ll)}$.

A tabloid $\tabloid{\TT} \in \Tabloidson{\mmu}$
is mapped to a tabloid by $\sigma\in S_{m}$ in $M^{\mmu}$.
Since $\widetilde I_{\mmu}(k;\Ll)$ is generated by
binomials $\tabloid{\TT} a^i_{\mmu,\Ll} - \zeta_l^{i k} \tabloid{\TT}$ of tabloids,
the representation matrix of $\sigma \in S_{m}$ for the basis
 $\overline{\Tabloidson{\mmu}}$ 
in $M^{\mmu}(k;\Ll)$
is
a matrix with entries in $\Set{1 , \zeta_l, \zeta_l^2, \ldots, \zeta_l^{l-1}}$
whose nonzero elements appear exactly once for each row and for
 each column.
Hence
\begin{gather*}
\charof (M^{\mmu}(k;\Ll))(\sigma)
=\sum_{j\in \ZZ/l\ZZ} \sum_{\TT\in\overline{\Tabloidson{\mmu}^{(\sigma,j;\Ll)}}}
 \zeta_l^{k j}
=\sum_{j\in \ZZ/l\ZZ}
 \numof{\overline{\Tabloidson{\mmu}^{(\sigma,j;\Ll)}}} \zeta_l^{k j}.
\end{gather*}
It follows from this equation that
\begin{align*}
\sum_{k \in \ZZ / l \ZZ} \zeta_l^{i k}\charof (M^{\mmu}(k;\Ll))(\sigma)
&=\sum_{k \in \ZZ / l \ZZ} \zeta_l^{i k}\sum_{j\in \ZZ/ l \ZZ}
 \numof{\overline{\Tabloidson{\mmu}^{(\sigma,j;\Ll)}}} \zeta_l^{k j}\\
&=\sum_{j\in \ZZ/ l\ZZ} \sum_{k \in \ZZ / l \ZZ}
\zeta_l^{k(i+j)} \numof{\overline{\Tabloidson{\mmu}^{(\sigma,j;\Ll)}}} .
\end{align*}
Since $ \sum_{k \in \ZZ / l \ZZ}\zeta_l^{kn} =0 $ for  $n \neq 0$,
this equation implies
\begin{gather*}
\sum_{k \in \ZZ / l \ZZ} \zeta_l^{i k}\charof (M^{\mmu}(k;\Ll))(\sigma)
= \sum_{k \in \ZZ / l \ZZ } 
 \numof{\overline{\Tabloidson{\mmu}^{(\sigma,-i;\Ll)}}}
= l\numof{\overline{\Tabloidson{\mmu}^{(\sigma,-i;\Ll)}}}.
\end{gather*}
Since 
$\numof{\Tabloidson{\mmu}^{(\sigma,-i)}}=\numof{\Braket{a_{\mmu,\Ll}}\overline{\Tabloidson{\mmu}^{(\sigma,-i;\Ll)}}}=l \numof{\overline{\Tabloidson{\mmu}^{(\sigma,-i;\Ll)}}}$,
we have 
\begin{gather*}
\sum_{k \in \ZZ / l \ZZ} \zeta_l^{i k}\charof (M^{\mmu}(k;\Ll))(\sigma)
=\numof{\Tabloidson{\mmu}^{(\sigma,-i;\Ll)}}.
\end{gather*}
\end{proof}

It follows from Lemma $\ref{equivlemma}$ and the equation $\eqref{eqofnum}$
 that Theorem \ref{gencharathm} and Theorem \ref{genidenthm} are equivalent.

We construct a bijection 
between 
marked $(\rho_\sigma,\ggamma,\Ll)$-tabloids on an $\Ll$-partition
$\mmu$ 
and
tabloids $\tabloid{\TT}$ on $\mmu$  satisfying 
$\sigma \tabloid{\TT} = \tabloid{\TT} a_{\mmu,\Ll}^{-1}$
to prove Theorem $\ref{genidenthm}$.

\begin{definition}
For a Young diagram $\rho\partitionof m$, we define   $n_{\rho,i}$,
 $N_{\rho,i}$, 
$\sigma_{\rho,i}$
 and
$\sigma_{\rho}$ by the following:
\begin{align*}
 n_{\rho,i}&=1+\sum_{j=1}^{i-1}\rho_j,\\
 N_{\rho,i}&=\Set{n_{\rho,i},n_{\rho,i}+1,\ldots,n_{\rho,i}+\rho_i-1 } \subset \Set{1,\ldots,m},\\
\sigma_{\rho,i}&=(n_{\rho,i},n_{\rho,i}+1,\ldots,n_{\rho,i}+\rho_i-1 )
 \in S_{m},\\
\sigma_{\rho}&=\sigma_{\rho,1}\sigma_{\rho,2}\sigma_{\rho,3}\cdots
 \in S_{m}.
\end{align*}
\end{definition}
For a Young diagram $\rho\partitionof m$,
it follows by definition that
$\bigcup_i N_{\rho,i} = \Set{1,\ldots, m}$,
$\numof{N_{\rho,i}}=\rho_i$
and the cycle type of  $\sigma_\rho$ is $\rho$.

\begin{definition}\label{defofphi2}
Let $\gamma_{h}$ be the order of $a_{\mu^{(h)},l_{h}}^{j}$.
For a marked $(\rho,\ggamma,\Ll)$-tabloid $(Y,c)$ on an $\Ll$-partition $\mmu$,
 $\tabloid{\varphi_{j}(Y,c)}$
denotes the tabloid obtained by the following:
\begin{itemize}
\item Put the number $n_{\rho,i}$ on a box in the $c(i)$-th row of 
the inverse image $Y^{-1}(\Set{i})$ for each $i$.
\item Put the number $\sigma_{\rho}n$ on a box in the $(\overline{c-j}+q l_{h})$-th row
of $\mu^{(h)}$ if the number $n$ is in the $(\overline{c}+q l_{h})$-th row of $\mu^{(h)}$,
where $\overline{c}$, $\overline{c-j} \in \ZZ/l_{h}\ZZ=\Set{1,\ldots,l_{h}}$ and $q\in \ZZ$.
\end{itemize}
We define $\varphi_j{(Y,c)}$ to be the numbering sorted in ascending order 
in each row of $\tabloid{\varphi_j{(Y,c)}}$.
\end{definition}

\begin{example}
For a marked  $((4,4,1),(2,1))$-tabloid
 $\left(\ \young(\boldII 2,22,11,\boldI 1)\ ,\ \young(\boldIII)\ \right)$,
\begin{gather*}
{\varphi_{1}\left(\ \young(\boldII 2,22,11,\boldI 1)\ ,\ 
\young(\boldIII)\ \right)}={\left(\ \young(57,68,24,13)\ ,\ \young(9)\ \right)}.
\end{gather*}
For a marked $((4,4,1),(4,1),(2,1))$-tabloid
$\left(\ \young(\boldII 2,11,22,\boldI1)\ ,\ \young(\boldIII)\ \right)$,
\begin{gather*}
{\varphi_2\left(\ \young(\boldII 2,11,22,\boldI1)\ ,\ \young(\boldIII)\
 \right)}={\left(\ \young(57,24,68,13)\ ,\ \young(9)\ \right)}.
\end{gather*}
\end{example}

Now we show that this correspondence $\varphi_j$ provides 
a realization of 
  Theorem \ref{genidenthm}.

\begin{lemma}\label{eigenlemma}
For a marked $(\rho,\ggamma,\Ll)$-tabloid $(Y,c)$ on 
an $\Ll$-partition  $\mmu$,
the tabloid $\tabloid{\varphi_{j}{(Y,c)}}$ satisfies
\begin{gather*}
\sigma_\rho\tabloid{\varphi_{j}{(Y,c)}}=\tabloid{\varphi_{j}{(Y,c)}}a_{\mmu,\Ll}^{-j},
\end{gather*}
where $\varphi_j$ is the one defined in Definition $\ref{defofphi2}$ and
 $\gamma_{h}$ is the order of $a_{\mu^{(h)},l_{h}}^{-j}$.
\end{lemma}
\begin{proof}
For the set $N_{\rho,i}$,
$\sigma_\rho$ acts as the cyclic permutation $\sigma_{\rho,i}$ of length
 $\rho_i$.
Since the set $(\varphi_j{(Y,c)})^{-1} (N_{\rho,i})$ of boxes
whose entries are in $N_{\rho,i}$
equals the inverse image $Y^{-1}(\Set{i})$,
we consider only $Y^{-1}(\Set{i})$ now.
Since we are allowed to change entries in the same row in $M^{\mmu}$,
it follows
from direct calculation that
$\sigma_\rho(\tabloid{\varphi_j{(Y,c)}})$ equals 
$\tabloid{\varphi_j{(Y,c+1)}}$
over $Y^{-1}(\Set{i})$,
where $(c+1)(i)=c(i)+1$ in $\ZZ/\gamma_{h}\ZZ$. 
Since $\sigma_\rho\tabloid{\varphi_j{(Y,c)}}$ equals 
$\tabloid{\varphi_j{(Y,c+1)}}$
over $Y^{-1}(\Set{i})$ for every $i$,
\begin{gather*}
\sigma_\rho\tabloid{\varphi_j{(Y,c)}}=\tabloid{\varphi_j{(Y,c+1)}}.
\end{gather*}
 
On the other hand, $a_{\mmu,\Ll}$ acts by
\begin{gather*}
\tabloid{\varphi_j{(Y,c)}} a_{\mmu,\Ll}^{j}=\tabloid{\varphi_j{(Y,c-1)}}.
\end{gather*}
Hence it follows that
\begin{gather*}
\sigma_\rho\tabloid{\varphi_j{(Y,c)}} 
=\tabloid{\varphi_j{(Y,c+1)}} 
=\tabloid{\varphi_j{(Y,c)}} a_{\mmu,\Ll}^{-j}.
\end{gather*}
\end{proof}

\begin{lemma}\label{shapelemma2}
Let a tabloid $\tabloid{\TT}$ on an $\Ll$-partition $\mmu$
satisfy 
$\sigma_\rho\tabloid{\TT}=\tabloid{\TT}a_{\mmu,\Ll}^{-j}$.
If $\TT^{-1}({n_{\rho,k}})$ is a box in the $(\overline{r}+l_{h}q)$-th 
row of $\mu^{(h)}$,
then  $n\in N_{\rho,k}$ is in the
 $(\overline{r-(n-n_{\rho,k})j}+l_{h}q)$-th  row of $\mu^{(h)}$,
where $\overline{r}$ and $\overline{r-(n-n_{\rho,k})j}\in \ZZ / l_{h}\ZZ =\Set{1,\ldots,l_{h}}$ and $q\in\ZZ$.
\end{lemma}

\begin{proof}
Let $\TT^{-1}({n_{\rho,k}})$ be a box in the $(\overline{r}+l_{h}q)$-th
 row in $\mu^{(h)}$,
 where $\overline{r} \in \ZZ/l_{h}\ZZ = \Set{1,2,\ldots,l_{h}}$.
Since
$\sigma_\rho\tabloid{\TT}=\tabloid{\TT} a_{\mmu,\Ll}^{-j},$
it follows, by the row where $n_{\rho,k}$ lies, that
$\sigma_\rho n_{\rho,k}$ is in the  $(\overline{r-j}+l_{h}q)$-th row  in $\mu^{(h)}$,
where $\overline{r-j} \in \ZZ/l_{h}\ZZ = \Set{1,2,\ldots,l_{h}}$.
Similarly,
we can   specify boxes where elements of
 $N_{\rho,k}$ lie.
Hence it follows that $n\in N_{\rho,k}$ lies in the $(\overline{r-(n-n_{\rho,k})j}+l_{h}q)$-th row,
where $\overline{r-(n-n_{\rho,k})j}\in \ZZ / l_{h}\ZZ =\Set{1,\ldots,l_{h}}$.
\end{proof}

\begin{lemma}\label{correspondlemma2}
If $\gamma_{h}$ is the order of $a_{\mu^{(h)},l_{h}}^{j}$,
our correspondence $\varphi_j$ provides
a bijection 
between 
marked $(\rho,\ggamma,\Ll)$-tabloids on an $\Ll$-partition $\mmu$
and
tabloids $\tabloid{\TT}$ on $\mmu$  satisfying 
$\sigma_\rho \tabloid{\TT} = \tabloid{\TT}a_{\mmu,\Ll}^{-j}$. 
\end{lemma}
\begin{proof}
We construct the inverse map of $\varphi_j$. 
Let   a tabloid $\tabloid{\TT}$ on $\mmu$  satisfy 
the equation 
$\sigma_\rho \tabloid{\TT} = \tabloid{\TT}a_{\mmu,\Ll}^{-j}$. 
Let $\TT$ be a numbering sorted in ascending order in each row
 of $\tabloid{\TT}$.

Let $\psi$ be a map from $\Set{1,\ldots,m}$
to $\NN$
which maps  $i \in N_{\rho,k}$
to $k$.
For $i<i'$, $n \in N_{\rho,i}$ is smaller than $n' \in N_{\rho,i'}$. 
Since  $\TT$ is sorted in ascending order in each row,
it follows from Lemma $\ref{shapelemma2}$ that
$\psi\circ \TT$ is 
a $(\rho,\ggamma,\Ll)$-tabloid.

Let $\psi_c (\TT)$ be
a map from $\Set{i | \rho_i \neq 0}$ to $\coprod_{h}\ZZ/l_{h}\ZZ$
 defined by
$\psi_c(i)=j \in \ZZ/\gamma_{h}\ZZ$
if $\TT^{-1}(n_{\rho,i})$ is in $\mu^{(h)}$ and
in the $j$-th row of the inverse image $\TT^{-1}(N_{(\rho,i)})$.
A pair $( \psi\circ \TT, \psi_c\TT)$
is  a marked $(\rho,\ggamma,\Ll)$-tabloid on $\mmu$.

Since $n_{\rho,i}$ is in the same row in  $\tabloid{\TT}$
as in $\tabloid{\varphi( \psi\circ \TT, \psi_c \TT)}$,
it follows from  Lemma $\ref{shapelemma2}$ that 
\begin{gather*}
\tabloid{\varphi( \psi\circ \TT, \psi_c (\TT))}=\tabloid{\TT}
\end{gather*}
and that
\begin{gather*}
(\psi\circ \varphi(Y,c),\psi_c (\varphi(Y,c)) )=(Y,c).
\end{gather*}
\end{proof}

Last we consider not only $\sigma_\rho$, but also general elements $\sigma$ whose
cycle type is $\rho$.
We explicitly give parameterizations of Theorem \ref{genidenthm} 
in the following theorem.

\begin{theorem}\label{lastlemma}
Let the  cycle type of $\sigma\in S_{m}$ be $\rho$ and let $\tau\in\ S_{m}$
satisfy $\tau \sigma_\rho \tau^{-1} =\sigma$.
Then 
the set
$\Set{\tabloid{\TT} \in \Tabloidson{\mmu}| 
\sigma \tabloid{\TT}=\tabloid{\TT} a_{\mmu,\Ll}^{-j}  }$ equals
\begin{gather*}
\Set{\tabloid{\tau\varphi_j(Y,c)}|
\text{$(Y,c)$ is a marked $(\rho,\ggamma,\Ll)$-tabloid on $\mmu$}}
\end{gather*}
for an $\Ll$-partition $\mmu$ of $\Mm$ and 
$\gamma_{h}$ is the order of $ a_{\mu^{(h)},l_{h}}^{j}$.
\end{theorem}
\begin{proof}
Since $\sigma = \tau \sigma_{\rho} \tau^{-1}$,
the equation
\begin{gather*}
\sigma\tabloid{\TT}=\tabloid{\TT} a_{\mmu,\Ll}^{-j}
\end{gather*}
is equivalent to the equation
\begin{gather*}
  \sigma_{\rho} \tau^{-1}\tabloid{\TT}=  \tau^{-1} \tabloid{\TT} a_{\mmu,\Ll}^{-j}.
\end{gather*}
From Lemma $\ref{correspondlemma2}$, there exists 
a marked  $(\rho,\ggamma,\Ll)$-tabloid $(Y,c)$ on $\mmu$
such that $\tabloid{\varphi_j(Y,c)} = \tau^{-1} \tabloid{\TT}$.
Hence this theorem follows.
\end{proof}

\begin{example}
Let $\mmu$,  $\Ll$ and $\rho$ be the same as ones in 
Example \ref{runningex}, i.e.,
$\mmu = ((2,2),4)$, $\Ll=(2,1)$ and $\rho=(4,2,2)$.
First we consider the case where $j=1$.
In this case,
\begin{gather*}
\left(\ \young(\boldI 1,11)\ ,\ \young(\boldII 2\boldIII 3)\ \right)
\end{gather*}
is a $(\rho,\Ll)$-tabloid on $\mmu$.
We have
\begin{gather*}
\varphi_1 \left(\ \young(\boldI 1,11)\ ,\ \young(\boldII 2\boldIII 3)\ \right)
=  \left(\ \young(13,24)\ ,\ \young(5678)\ \right).
\end{gather*}
Since $\sigma_\rho=(1,2,3,4)(5,6)(7,8)$ acts as
\begin{gather*}
(1234)(56)(78) \left(\ \young(13,24)\ ,\ \young(5678)\ \right)
= \left(\ \young(24,31)\ ,\ \young(6587)\ \right)
\end{gather*}
and 
\begin{align*}
 \tabloid{\left(\ \young(24,31)\ ,\ \young(6587)\ \right)}
&=\tabloid{\left(\ \young(24,13)\ ,\ \young(5678)\ \right)}\\
&=\tabloid{\left(\ \young(13,24)\ ,\ \young(5678)\ \right)} a_{\mmu,\Ll}^{-1},
\end{align*}
it follows that
\begin{gather*}
\sigma_\rho \tabloid{\varphi_1 \left(\ \young(\boldI 1,11)\ ,\
 \young(\boldII 2\boldIII 3)\ \right)}
= \tabloid{\varphi_1 \left(\ \young(\boldI 1,11)\ ,
\ \young(\boldII 2\boldIII 3)\ \right)} a_{\mmu,\Ll}^{-1}.
\end{gather*}
Next we consider the case where $j=2$.
In this case, 
\begin{gather*}
\left(\ \young(\boldII 2,\boldIII 3)\ ,\ \young(\boldI 111)\ \right)
\end{gather*}
is a $(\rho,(1,1),\Ll)$-tabloid on $\mmu$.
We have
\begin{gather*}
\varphi_2 \left(\ \young(\boldII 2,\boldIII 3)\ ,\ \young(\boldI 111)\ \right)
=  \left(\ \young(56,78)\ ,\ \young(1234)\ \right).
\end{gather*}
Since $\sigma_\rho$ acts as
\begin{align*}
(1234)(56)(78) \tabloid{\left(\ \young(56,78)\ ,\ \young(1234)\ \right)}
&= \tabloid{ \left(\ \young(65,87)\ ,\ \young(2341)\ \right)}\\
&= \tabloid{ \left(\ \young(56,78)\ ,\ \young(1234)\ \right)}\\
\end{align*}
and $ a_{\mmu,\Ll}^{-2} = \varepsilon \in S_{8}$,
it follows that
\begin{gather*}
\sigma_\rho \tabloid{\varphi_1 \left(\ \young(\boldII 2,\boldIII 3)\ ,\
 \young(\boldI 111)\ \right)}
= \tabloid{\varphi_1 \left(\ \young(\boldII 2,\boldIII 3)\ ,
\ \young(\boldI 111)\ \right)} a_{\mmu,\Ll}^{-2}.
\end{gather*}

\end{example}

\end{document}